\documentclass[12pt]{article}
\usepackage[cp1251]{inputenc}
\usepackage[russian]{babel}
\usepackage{amsfonts,amsmath,amssymb,amsthm,graphicx,afterpage}
\usepackage{epsfig}
\graphicspath{{morefig/}{figures/}{figures00/}{figures02/}{figures05/}{figlms/}}

\DeclareGraphicsExtensions{.pdf,.png,.jpg}

\def\a{\alpha}
\def\s{\varepsilon}
\def\b{\beta}
\def\c{\gamma}
\def\cc{\hat{\gamma}}
\def\TrivSetKappa{(\kappa,\xi(\kappa))}

\def\o{\kappa}
\def\Triv{\overline{\mathcal{K}}}
\def\TrivS{\overline{\mathcal{K}'}}
\def\TrivSet{\mathcal{K}}
\def\TrivSetS{\mathcal{K}'}
\def\B{K}

\def\R{{\mathbb R}}

\long\def\comment#1\endcomment{}

\newtheoremstyle{mydefinition}
{3pt}
{3pt}
{\normalfont}
{\parindent}
{\bfseries}
{.}
{ }
{}

\theoremstyle{plain}
\newtheorem{theorem}{Теорема}

\newtheorem{lemma}[theorem]{Лемма}

\newtheorem{approvalo}[theorem]{Утверждение}

	\title{Элементарное доказательство существования полинома Конвея \thanks{Автор поддержен грантом РФФИ  19-01-00169.}}

	\author{ Гараев Тимур Рустемович}
	
	\date{}

\begin{document}	
	\maketitle
	

\par

\begin{center}
\section{Формулировка теоремы}	
\end{center}

Назовем \textit{неупорядоченной диаграммой\footnote{Определение \textit{диаграммы} приведено в \cite{TurGor}.}} плоскую диаграмму ориентированного зацепления с неупорядоченными компонентами.

\begin{figure}[h]
	\begin{center}
	\includegraphics[scale=0.7]{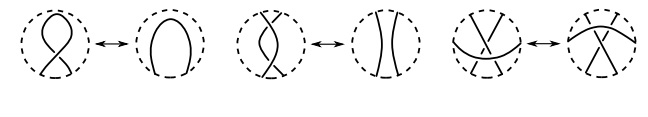}
	
	 $R_1$ \ \ \ \ \ \ \ \ \ \ \ \ \ \ \ \ \ \ \ \ \ \ $R_2$ \ \ \ \ \ \ \ \ \ \ \ \ \ \ \ \ \ \ \ \ \ \ $R_3$
	 	\end{center}
	 \caption{}
\label{ris:conversion}
\end{figure}

Обозначим через $R_1$, $R_2$ и $R_3$ движения Рейдемейстера неупорядоченных диаграмм, изображенные на рис. \ref{ris:conversion}. Обозначим через $R_0$ плоскую изотопию неупорядоченных диаграмм.

Будем называть неупорядоченные диаграммы $\B$ и $\B'$ \textit{эквивалентными}, если $\B'$ получена из  $\B$ применением движений Рейдемейстера и плоских изотопий.

Будем называть неупорядоченную диаграмму \textit{тривиальной}, если она эквивалентна неупорядоченной диаграмме, состоящей из попарно непересекающихся несамопересекающихся ломаных. 

Будем называть функцию  $f$ на множестве неупорядоченных диаграмм \textit{инвариантом неупорядоченных диаграмм}, если $f$ совпадает на эквивалентных неупорядоченных диаграммах.

\begin{figure}[h]
	\begin{center}
		\includegraphics[scale=0.7]{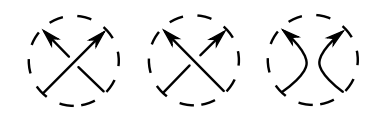}

		$K_+$ \ \ \ \ \ \ \ \ \ \ \ \ $ K_-$ \ \ \ \ \ \ \ \ \ \ \ \ $K_0$
		
		
	\end{center}
	\caption{}
	\label{ris:K}
\end{figure}

\begin{theorem}\label{p:Theor}
	\textit{Существует бесконечная последовательность $c_{-1}=0,c_0,c_1,\ldots$ инвариантов неупорядоченных диаграмм, принимающая значение $c_0=1, c_1=0, c_2=0,\ldots$ на неупорядоченной тривиальной однокомпонентной диаграмме, и такая, что для любых трех неупорядоченных диаграмм $K_+$, $K_-$ и $K_0$, отличающиеся внутри пунктирных кругов, как показано на рис. 2, и совпадающие вне этих кругов и для любого $n\ge0$ выполняется равенство
		$$c_n(K_+)-c_n(K_-)=c_{n-1}(K_0).$$}
\end{theorem}

Известные доказательства  теоремы \ref{p:Theor} приведены, например, в \cite{Ale30}, \cite{Con}, \cite{Ma18} и \cite{Kauffman}.
Приводимое здесь доказательство не является более коротким, но оно основано на другой идее (см. $\S 2$ и $\S 3$), которая может быть применена к элементарному построению полинома HOMFLY или, возможно, к элементарному построению инвариантов Васильева-Гусарова.

Автор благодарен А. Скопенкову и Е. Морозову за помощь в составлении текста.

\begin{center}
	\section{Доказательство теоремы \ref{p:Theor} с использованием леммы \ref{p:GreLem}}
\end{center}

Будем называть \textit{упорядоченной диаграммой} плоскую диаграмму ориентированного зацепления с упорядоченными компонентами.

\begin{figure}[h]
	\begin{center}

		\includegraphics[scale=1.5]{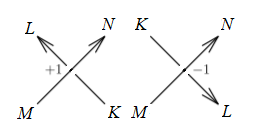}
		
	\end{center}
	\caption{}
	\label{ris:arrow}
\end{figure}

Пусть перекресток $a$ неупорядоченной (упорядоченной) диаграммы $\B$ образован пересечением отрезков $MN$ и $KL$. Пусть $(\overrightarrow{MN}, \overrightarrow{KL})$ — упорядоченная пара пересекающихся векторов, где $MN$ проходит над $KL$. Определим \textit{знак перекрестка} $a$ как $+1$, если $MNL$ ориентирован против часовой
стрелки и как $-1$ иначе (рис. \ref{ris:arrow}). 

Будем называть перекресток $a$ со знаком  $+1$ \textit{положительным}.  Будем называть перекресток $a$ со знаком  $-1$ \textit{отрицательным}. 
Обозначим знак перекрестка $a$  неупорядоченной (упорядоченной) диаграммы $\B$ через $\s(a,\B)$. 

Предположим, что перекресток $a$ неупорядоченной (упорядоченной) диаграммы $\B$ образован пересечением отрезков $MN$ и $KL$, где $KL$ лежит под $MN$.
Назовем \textit{изменением перекрестка} $a$ замену неупорядоченной (упорядоченной) диаграммы $\B$ на неупорядоченную (упорядоченную) диаграмму  $\B_a$, у которой отрезок $MN$ лежит под отрезком $KL$ и которая отличается от $\B$ только в перекрестке $a$.

\begin{figure}[h]
	\begin{center}

		\includegraphics[scale=1]{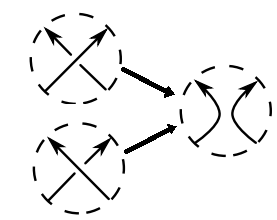}
		
	\end{center}
	\caption{}
	\label{ris:delete}
\end{figure}

Назовем \textit{удалением перекрестка} операцию, переводящую неупорядоченную (упорядоченную) диаграмму в неупорядоченную диаграмму, как на рис. \ref{ris:delete}.

Обозначим через $\B_{[a]}$ неупорядоченную диаграмму, полученную из неупорядоченной (упорядоченной) диаграммы $\B$ удалением перекрестка $a$.


Обозначим через $\B_a$ неупорядоченную (упорядоченную) диаграмму, полученную из неупорядоченной (упорядоченной) диаграммы $\B$ изменением перекрестка $a$.

\smallskip

\begin{lemma}\label{p:MGreLem}
	 Пусть дана такая пара $(\a,\b)$ инвариантов неупорядоченных диаграмм, что $\b(\B)=0$ для любой тривиальной неупорядоченной диаграммы $\B$, у которой больше одной компоненты и что $$\b(\B)-\b(\B_a)=\s(a,\B)\a(\B_{[a]})$$ для любой неупорядоченной диаграммы $\B$. Существует такой инвариант $\gamma$ неупорядоченных диаграмм, что
	  
	$\bullet$ $\c(\B)-\c(\B_a)=\s(a,\B)\a(\B_{[a]})$ для любой неупорядоченной диаграммы $\B$;
	
	$\bullet$ $\gamma(\B)=0$ для любой тривиальной неупорядоченной диаграммы $\B$.
\end{lemma}
	 
\textit{Доказательство теоремы \ref{p:Theor} с использованием леммы \ref{p:MGreLem}.}
Определим $c_0(\B)=1$ для любой неупорядоченной диаграммы $\B$, у которой одна компонента, и $c_0(\B')=0$ для любой неупорядоченной диаграммы $\B'$, у которой больше одной компоненты. Тогда пара $(c_{-1},c_{0})$ является конвеевской. Пусть существует последовательность $c_{-1},c_0,\ldots,c_n$ инвариантов неупорядоченных диаграмм, удовлетворяющая условиям теоремы \ref{p:Theor}. Тогда пара $(c_{n-1},c_n)$ является конвеевской и $c_n(\B)=0$ для любой неупорядоченной тривиальной диаграммы $\B$, у которой больше одной компоненты. По лемме \ref{p:MGreLem} существует инвариант $c_{n+1}$ неупорядоченных диаграмм такой, что пара $(c_n,c_{n+1})$ является конвеевской и $c_{n+1}(\B)=0$ для любой неупорядоченной тривиальной диаграммы $\B$. 
Так строится бесконечная последовательность, удовлетворяющая условию теоремы \ref{p:Theor}.  $\Box$ 

\smallskip




 

Движение Рейдемейстера, эквивалентность, тривиальность и инвариантность для упорядоченных диаграмм определяются по аналогии с неупорядоченными диаграммами.

Будем называть пару $(\a,\b)$ вещественнозначных функций, где функция $\a$ задана на множестве неупорядоченных диаграмм, а функция $\b$ задана на множестве упорядоченных диаграмм, \textit{конвеевской}, если
для любой упорядоченной диаграммы $\B$ выполняется $\b(\B_{a})-\b(\B)=\s(a,\B_{a})\a(\B_{[a]})$.
 
 \begin{lemma}\label{p:GreLem}
 	
Дана конвеевская пара $(\a,\b)$ инвариантов, где $\b(\B)=0$  для любой тривиальной неупорядоченной диаграммы $\B$, у которой больше одной компоненты. Существует такой инвариант $\c$ упорядоченных диаграмм, что
 	
$\bullet$ пара $(\b,\c)$ является конвеевской;
 	
$\bullet$ $\c(\B)=0$ для любой тривиальной упорядоченной диаграммы $\B$.
 \end{lemma}

\textit{Доказательство леммы \ref{p:MGreLem} с использованием леммы \ref{p:GreLem}.} 
Для доказательства леммы \ref{p:MGreLem} достаточно показать, что инвариант $\c$ упорядоченных диаграмм из леммы \ref{p:GreLem} совпадает на диаграммах, отличающихся только порядком компонент.

Обозначим через $\B$ упорядоченную диаграмму с двумя или более компонентами. По теореме 3.8 из \cite{SP}, из упорядоченной диаграммы $\B$ можно получить, применяя операции изменения перекрестков, упорядоченную тривиальную диаграмму. 
Поэтому можно обозначить через $u(\B)$ минимальное число изменений перекрестков, необходимое для получения упорядоченной тривиальной диаграммы из упорядоченной диаграммы $\B$.
Докажем лемму \ref{p:MGreLem} индукцией 
по $u(\B)$. Если $u(\B) = 0$,
то $\B$ — упорядоченная тривиальная диаграмма. Обозначим через $\B'$ упорядоченную диаграмму, полученную из $\B$ изменением порядка компонент. Тогда $\c(\B)=\c(\B')=0$. Предположим, что $u(\B) > 0$. Обозначим через $a$ такой перекресток диаграммы $\B$, что $u(\B_a) < u(\B)$. Тогда
$$\c(\B) - \c(\B_a) = \s(a,\B)\b(\B_{[a]}) \quad\text{и}\quad \c(\B') - \c(\B_a') = \s(a,\B')\b(\B_{[a]}').$$
Заметим, что упорядоченные диаграммы $\B_a$ и $\B_a'$
совпадают с точностью до порядка компонент, а
неупорядоченные диаграммы $\B_{[a]}$ и $\B_{[a]}'$ равны. Следовательно, $\b(\B_{[a]})=\b(\B_{[a]}')$. Так как $u(\B_a) < u(\B)$, по предположению индукции 
имеем $\c(\B_a) = \c(\B'_a)$. Следовательно, $\c(\B) = \c(\B').$ $\Box$

\smallskip
 
 \begin{center}
 	\section{Построение функции $\c$ из леммы \ref{p:GreLem}}
 \end{center}

Будем называть упорядоченную диаграмму \textit{диаграммой}.
	
Будем называть \textit{контуром} замкнутую ориентированную ломаную на плоскости. 

Будем говорить, что точки на плоскости находятся \textit{в общем положении}, если никакие три точки не лежат на одной прямой и если никакие три прямые, проходящие через шесть разных точек, не пересекаются в одной точке. 


Будем называть упорядоченный набор контуров, множество вершин которых находится в общем положении, \textit{набором контуров}.

Будем называть пару $(\kappa,\xi)$, где $\o=(\o_1,\ldots,\o_n)$ --- набор контуров, а $\xi=(\xi_1,\ldots,\xi_n)$ --- такой упорядоченный набор точек, что $\xi_i \in \o_i$ и $\xi_i$ не является точкой самопересечения контура $\o_i$ для $i=1,\ldots,n$, \textit{набором контуров с отмеченными точками}.

Будем говорить, что компонента $\bar{\o_i}$ диаграммы $\B$ \textit{лежит под} компонентой $\bar{\o_j}$ той же диаграммы, если любой отрезок компоненты $\bar{\o_i}$ либо не пересекается с компонентой $\bar{\o_j}$, либо лежит под отрезком, лежащим в компоненте $\bar{\o_j}$.
 
Для набора $\TrivSet=(\kappa,\xi)$ контуров с отмеченными точками обозначим через $\Triv$ диаграмму, получающуюся  
из набора $\o$ добавлением информации, какой из пересекающихся отрезков лежит ниже, а какой выше,  с компонентами $\bar \o_1,\ldots,\bar \o_n$ такую, что
	  
$\bullet$ для любого $i<j$, компонента $\bar{\o_i}$ диаграммы $\Triv$ лежит под компонентой $\bar{\o_j}$ той же диаграммы;
	  
$\bullet$ для любого $i$, при обходе компоненты $\bar{\o_i}$, начиная с точки $\xi_i$, в соответствии c ориентацией компоненты $\bar{\o_i}$, каждый отрезок $KL$ компоненты $\bar{\o_i}$ лежит выше любого отрезка, который был пройден раньше отрезка $KL$. 

Нетрудно видеть, что диаграмма $\Triv$ тривиальная\footnote{Построение диаграммы $\Triv$ схоже с описанным в доказательстве теоремы 3.8 из \cite{SP}.}.
	  
Обозначим через $x_1,\ldots, x_j$ некоторые перекрестки (не обязательно различные) некоторой диаграммы $\B$. Обозначим $x:=(x_1,x_2, \ldots,x_j)$.
Обозначим через $\B_{x}$ диаграмму, полученную из диаграммы $\B$ изменением перекрестков $x_1,\ldots,x_j$.


Каждый набор $\kappa$ контуров 
произвольно дополним до набора $\TrivSetKappa$ контуров с отмеченными точками.
Обозначим через $\mathbb{B}$ множество всех пар $(x,\kappa)$, где $\kappa$ --- набор контуров, а $x$ ---  упорядоченный набор 
из некоторых перекрестков (не обязательно различных) диаграммы $\overline{\TrivSetKappa}$. Если набор $x$ не содержит элементов, то пару $(x, \kappa)$ будем записывать как $(\kappa)$.
	    
Зададим функцию $\cc:\mathbb{B} \rightarrow \R$ следующим образом: 
$$\cc(\kappa):=0\quad\text{и}\quad\cc(x,a,\kappa):= \cc(x,\kappa)-\s(a,\overline{\TrivSetKappa}_{x})\b(\overline{\TrivSetKappa}_{x,[a]}).$$	

Для каждой диаграммы $\B$, которая получается из упорядоченного набора $\kappa$ контуров добавлением информации, какой из пересекающихся отрезков проходит ниже, а какой выше, выберем такой упорядоченный набор $y=y(\B,\overline{\TrivSetKappa})$ перекрестков, что $\overline{\TrivSetKappa}_y=\B$.
Зададим функцию $\c$ из леммы \ref{p:GreLem} формулой $$\c(\B):=\cc(y,\kappa).$$ Значение $\c(\B)$ зависит от выбора набора $\TrivSetKappa$ по $\kappa$, но этот выбор считается фиксированным.
 Функция $\c$ корректно определена, а именно $\c(\B)$ не зависит от выбора $y(\B,\overline{\TrivSetKappa})$, ввиду утверждения \ref{p:Kor}.

	


\smallskip
	
\begin{approvalo}\label{p:Kor}
	\textit{
	Пусть $\kappa$ --- набор контуров. Если наборы $x,y$ перекрестков диаграммы $\overline{\TrivSetKappa}$ такие, что $\overline{\TrivSetKappa}_x=\overline{\TrivSetKappa}_{y}$, то $\cc(x,\kappa)=\cc(y,\kappa)$. 
	}
\end{approvalo}

\textit{Доказательство.}
Обозначим $\TrivSet=\TrivSetKappa$.
Обозначим через $z_1,\ldots, z_j,a,b$ перекрестки диаграммы $\Triv$, где $a \neq b$ (в остальном перекрестки не обязательно различны). Обозначим $z:=(z_1,z_2, \ldots,z_j)$.
 
Имеем
$$\cc(z,a,a,\kappa)=\cc(z,a,\kappa)-\s(a,\Triv_{z,a})\b(\Triv_{z,a,[a]})=$$
$$=\cc(z,\kappa)-\s(a,\Triv_{z})\b(\Triv_{z,[a]})-\s(a,\Triv_{z,a})\b(\Triv_{z,a,[a]})=\cc(z,\kappa),$$
Так как $\s(a,\Triv_{z})=-\s(a,\Triv_{z,a})$. Имеем
$$\cc(z,a,b,\kappa)=\cc(z,a,\kappa)-\s(b,\Triv_{z,a})\b(\Triv_{z,a,[b]})=$$
$$=\cc(z,\kappa)-\s(a,\Triv_{z})\b(\Triv_{z,[a]})-\s(b,\Triv_{z,a})\b(\Triv_{z,a,[b]})=$$
$$=\cc(z,\kappa)-\s(a,\Triv_{z})\b(\Triv_{z,[a]})-\s(b,\Triv_{z,a})(\b(\Triv_{z,[b]})-\s(a,\Triv_{z,[b]})\a(\Triv_{z,[a],[b]}))=$$
$$=\cc(z,\kappa)-\s(a,\Triv_{z})\b(\Triv_{z,[a]})-\s(b,\Triv_{z})\b(\Triv_{z,[b]})+\s(b,\Triv_{z})\s(a,\Triv_{z})\a(\Triv_{z,[a],[b]}).$$
Так как $\Triv_{z,[a],[b]}=\Triv_{z,[b],[a]}$, то последнее выражение симметрично по $a$ и $b$. Следовательно, 
$$\cc(z,a,b,\kappa)=\cc(z,b,a,\kappa).$$    
Поэтому $\cc(z,\kappa)$ не зависит от порядка, в котором идут элементы в упорядоченном наборе $z$. Из этого, а также из равенства $\cc(z,a,a,\kappa)=\cc(z,\kappa)$ следует, что $\cc(z,\kappa)$ зависит только от четности вхождения каждого перекрестка в набор $z$. 

Так как $\Triv_x=\Triv_y$, то знак перекрестка $a$ диаграммы $\Triv_x$ совпадает со знаком  перекрестка $a$ диаграммы $\Triv_{y}$. 
Поэтому четность вхождения перекрестка $a$ в $x$ совпадает с четностью вхождения перекрестка $a$ в $y$. 
Следовательно, $\cc(x,\Triv)=\cc(y,\Triv).$
$\Box$

\smallskip

\textit{Доказательство конвеевости пары $(\b,\c)$.}
Нетрудно видеть, что существует такой набор $\kappa$ контуров, что диаграммы $\B$ и $\B_a$ получаются из $\kappa$ добавлением информации, какой из пересекающихся отрезков проходит ниже, а какой выше. 
Выберем такой набор $x$, что
$\overline{\TrivSetKappa}_{x}=\B$. Тогда $\overline{\TrivSetKappa}_{x,a}=\B_a$ и $\overline{\TrivSetKappa}_{x,[a]}=\B_{[a]}$. Имеем $$\c(\B)-\c(\B_a)=\cc(x,\kappa)-\cc(x,a,\kappa)=\s(a,\overline{\TrivSetKappa}_{x})\b(\overline{\TrivSetKappa}_{x,[a]})=\s(a,\B)\b(\B_{[a]}).$$ Следовательно, пара функций $(\b,\c)$ является конвеевской. $\Box$

\smallskip
	   
Для доказательства леммы \ref{p:GreLem} теперь достаточно показать инвариантность функции $\c$, что будет сделано в $\S 4$. 
	

\begin{center}
	\section{Доказательство инвариантности функции $\c$}	 
\end{center}

\smallskip 

Будем называть две диаграммы $\B$ и $\B'$ \textit{схожими}, если $\B'$ может быть получена из $\B$  операциями изменения перекрестка. 

Назовем  \textit{дугой диаграммы} $\B$ ломаную, принадлежащую диаграмме $\B$, не содержащую перекрестков и ограниченную двумя перекрестками. Обозначим дугу диаграммы, ограниченную перекрестками $a$ и $b$, через $ab$. 

\smallskip 

\begin{approvalo}\label{p:hardest}
	\textit{
		Для любого набора $\TrivSet=(\kappa,\xi)$ контуров с отмеченными точками выполняется $\c(\Triv)=0$. 
	}
\end{approvalo}

\smallskip 

\textit{Доказательство утверждения \ref{p:hardest}.}
Из построения функции $\c$ следует, что для диаграммы $\Triv$ существует набор $\mathcal{K}'=(\kappa,(\xi_1',\ldots,\xi_n'))$ контуров с отмеченными точками для которого $\c(\overline{\mathcal{K}'})=0$. 
Для доказательства утверждения \ref{p:hardest} достаточно рассмотреть случай, когда набор $\TrivSet$ получен из $\mathcal{K}'$ заменой отмеченной точки $\xi_i'$ на другую точку $\xi_i$.

Если точки $\xi_i'$ и $\xi_i$ принадлежат одной и той же дуге, то $\overline{\mathcal{K}'}=\Triv$.

\begin{figure}[h]
	\begin{center}

		\includegraphics[scale=1.5]{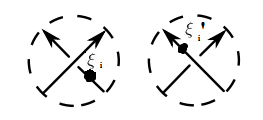}
		
		$\Triv$ \ \ \ \ \ \ \ \ \ \ \ \ \ \ \ \ \ \ \ \ \ \ \ \ \ \ \ \ \ \ \ \ $\TrivS$
	\end{center}
\end{figure}

Если точки $\xi_i'$ и $\xi_i$ принадлежат разным дугам, то достаточно рассмотреть случай, когда точка   $\xi_i'$ принадлежит дуге $ab$, а точка $\xi_i$ принадлежит дуге $bc$. Если перекресток $b$ образован пересечением двух разных компонент, то $\overline{\mathcal{K}'}=\Triv$. Если перекресток $b$ образован самопересечением компоненты $\bar\o_i'$, то  $$\s(b,\overline{\mathcal{K}'})\b(\overline{\mathcal{K}'}_{[b]})=\c(\overline{\mathcal{K}'})-\c(\overline{\mathcal{K}'}_b)=\c(\overline{\mathcal{K}'})-\c(\Triv).$$
Нетрудно видеть, что неупорядоченная диаграмма $\TrivS_{[b]}$ является неупорядоченной тривиальной диаграммой, у которой больше одной компоненты. Следовательно, $\b(\TrivS_{[b]})=0$. Следовательно, $\c(\overline{\mathcal{K}'})=\c(\overline{\mathcal{K}'}_b)=\c(\Triv)=0.$ 
$\Box$ 

\smallskip

Обозначим через $\B_{R_i}$ любую диаграмму (неупорядоченную диаграмму), полученную из диаграммы (неупорядоченной диаграммы) $\B$ преобразованием $R_i$ для $i = 0,1,2,3$.

\begin{figure}[h]
	\begin{center}

		\includegraphics[scale=0.4]{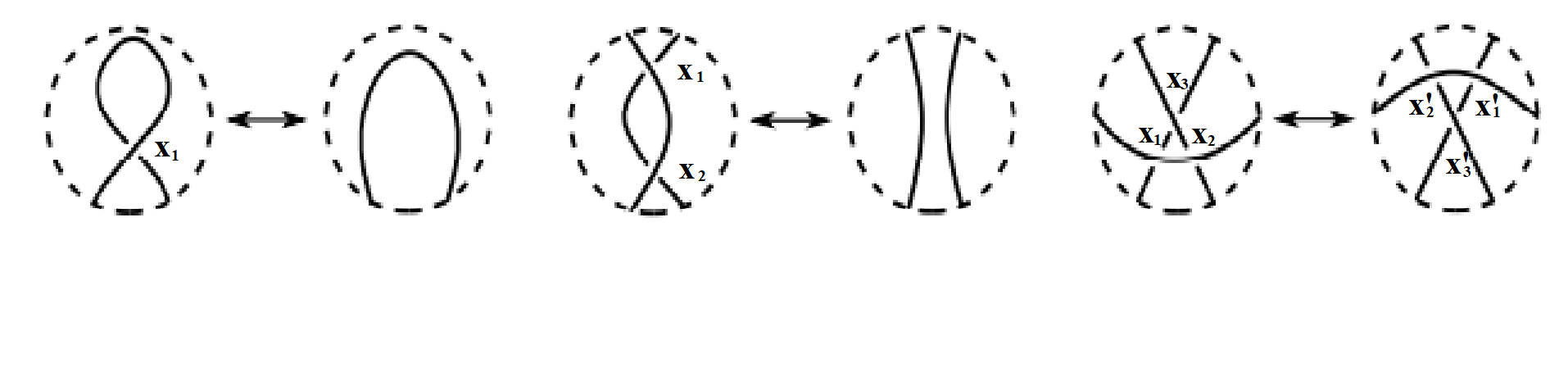}
		
	\end{center}
	\caption{}
	\label{ris:Perecret}
\end{figure}

Будем говорить, что перекрестки $x_1$ и $x_2$ с рис. \ref{ris:Perecret} \textit{участвуют} в движении Рейдемейстера $R_i$ для $i=1,2,3$.
Если перекрестки $x_1,\ldots,x_j$ (не обязательно различные) диаграммы (неупорядоченной диаграммы) $\B$ не участвуют в движении Рейдемейстера $R_i$ для $i=1,2,3$, то будем
обозначим через $x_1',\ldots,x_j'$ перекрестки диаграммы (неупорядоченной диаграммы) $\B_{R_i}$, в которые переходят перекрестки $x_1,\ldots,x_j$ при $R_i$.
Обозначим через $x_1',\ldots,x_j'$ перекрестки диаграммы (неупорядоченной диаграммы) $\B_{R_0}$, в которые переходят перекрестки $x_1,\ldots,x_j$ (не обязательно различные) при плоской изотопии $R_0$. Будем говорить, что перекрестки $x_1,\ldots,x_j$ не участвуют в плоской изотопии.  Обозначим $x':=(x_1',\ldots,x_j')$.

\begin{approvalo}\label{p:Mov}
	\textit{
	Если $\c(\B)=\c(\B_{R_i})$, где $i = 0,1,2,3$ и перекресток  $a$ не участвует в преобразовании $R_i$, то $\c(\B_{a})=\c(\B_{R_i,a'})$.
	}
\end{approvalo}

\textit{Доказательство.}	Имеем $$\c(\B_{a})=\c(\B)+\s(a,\B_{a})\b(\B_{[a]})=\c(\B_{R_i})+\s(a',\B_{R_i,a'})\b(\B_{R_i,[a']})=\c(\B_{R_i,a'}).$$ Так как $\c(\B)=\c(\B_{R_i})$ и $\b(\B_{[a]})=\b(\B_{R_i,[a']})$, то $\c(\B_{a})=\c(\B_{R_i,a'}). \Box$

\smallskip

Покажем инвариантность функции $\c$ относительно каждого из преобразований $R_i$, где $i=0,1,2,3$.

\smallskip

\textit{Доказательство равенства $\c(\B)=\c(\B_{R_0})$.} 
Выберем такой набор $\TrivSet=(\kappa,\xi)$ контуров с отмеченными точками, что диаграмма $\Triv$ схожа с диаграммой $\B$. Обозначим через $\TrivSetS=(\kappa',\xi')$ такой набор контуров с отмеченными точками, что диаграмма $\TrivS$ схожа с диаграммой $\B_{R_0}$ и $\TrivS$ получена из $\Triv$ преобразованием $R_0$.
По утверждению \ref{p:hardest}   $0=\c(\Triv)=\c(\TrivS)=\c(\Triv_{R_0})$. Тогда по утверждению \ref{p:Mov}  $\c(\Triv_x)=\c(\Triv_{R_0,x'})=\c(\B)=\c(\B_{R_0}).$ $\Box$

\smallskip

	
		
\textit{Доказательство равенства $\c(\B)=\c(\B_{R_1})$.}
Выберем такой набор $\TrivSet=(\kappa,\xi)$ контуров с отмеченными точками, что диаграмма $\Triv$ схожа с диаграммой $\B$. 
Обозначим через $\TrivSetS=(\kappa',\xi')$ такой набор контуров с отмеченными точками, что диаграмма $\TrivS$ схожа с диаграммой $\B_{R_1}$ и $\TrivS$ получена из $\Triv$ движением $R_1$. Обозначим через $a$ перекресток, участвующий при $R_1$, диаграммы $\B_{R_1}$. Нетрудно видеть, что диаграмма $\TrivS_{[a]}$ тривиальная и имеет более одной компоненты. Следовательно, $\b(\TrivS_{[a]})=0$. Имеем $$ \c(\TrivS)-\c(\TrivS_a)=\s(a,\TrivS)\b(\TrivS_{[a]})=0.$$ Следовательно, $\c(\TrivS)=\c(\TrivS_{a})$. Заметим, что либо у диаграммы $\TrivS_{a}$, либо у диаграммы $\TrivS$ перекресток, участвующий при $R_1$, совпадают с перекрестком у диаграммы $\B_{R_1}$, участвующим при $R_1$. Обозначим соответствующую диаграмму через $\Triv_{R_1}$.

			
Выберем такое $x$, что $\B=\Triv_x$ и $x$ не содержит перекрестка $a$. По утверждению \ref{p:hardest}  $0=\c(\Triv)=\c(\TrivS)=\c(\Triv_{R_1})$. Тогда по утверждению \ref{p:Mov}  $\c(\Triv_x)=\c(\Triv_{R_1,x'})=\c(\B)=\c(\B_{R_1})$. $\Box$

\smallskip
		
\begin{figure}[h]
	\begin{center}

			\ \ \ \ \ \ \includegraphics[scale=0.4]{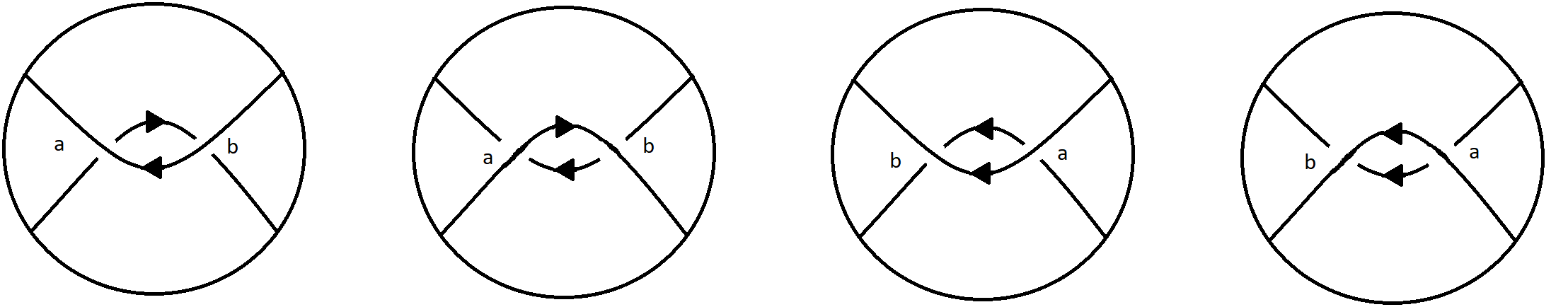}
	%
	\end{center}
		\caption{}
	\label{ris:R_2}
\end{figure}

		
\textit{Доказательство равенства $\c(\B)=\c(\B_{R_2})$.}
Выберем такой набор $\TrivSet=(\kappa,\xi)$ контуров с отмеченными точками, что диаграмма $\Triv$ схожа с диаграммой $\B$. 
Обозначим через $\TrivSetS=(\kappa',\xi')$ такой набор контуров с отмеченными точками, что диаграмма $\TrivS$ схожа с диаграммой $\B_{R_2}$ и $\TrivS$ получена из $\Triv$ движением $R_2$. Обозначим через $a$ и $b$ перекрестки, участвующие при $R_2$, диаграммы $\B_{R_2}$ как на рис. \ref{ris:R_2}. Имеем
		
$$\c(\TrivS)=\c(\TrivS_{b})+\s(b,\TrivS)\b(\TrivS_{[b]})=\c(\TrivS_{a,b,a})+\s(a,\TrivS_{a,b})\b(\TrivS_{a,b,[a]})=\c(\TrivS_{a,b}).$$

Следовательно, $\c(\TrivS)=\c(\TrivS_{a,b})$. Заметим, что либо у диаграммы $\TrivS_{a,b}$, либо у диаграммы $\TrivS$ перекрестки, участвующие при $R_2$, совпадают с перекрестками у диаграммы $\B_{R_2}$, участвующими при $R_2$. Обозначим соответствующую диаграмму через $\Triv_{R_2}$. 

 Выберем такой $x$, что $\Triv_x=\B$ и $x$ не содержит перекрестков $a$ и $b$. По утверждению \ref{p:hardest}  $0=\c(\Triv)=\c(\TrivS)=\c(\Triv_{R_2})$.
Тогда по утверждению \ref{p:Mov}  $\c(\Triv_x)=\c(\Triv_{R_2,x'})=\c(\B)=\c(\B_{R_2})$. $\Box$

\smallskip
		
\textit{Доказательство равенства $\c(\B)=\c(\B_{R_3})$.}
Обозначим через $x_1,x_2,x_3$ перекрестки диаграммы $\B$, участвующие в преобразовании $R_3$ как на рис. 5.

Построим такой набор $\TrivSet=(\kappa,\xi)$ контуров с отмеченными точками, что диаграмма $\Triv$ схожа с диаграммой $\B$ и отмеченные точки не лежат на дугах $x_1x_2$, $x_2x_3$, $x_3x_1$ диаграммы $\Triv$. Обозначим через $\TrivSetS=(\kappa',\xi')$ такой набор контуров с отмеченными точками, что диаграмма $\TrivS$ схожа с диаграммой $\B_{R_3}$ и $\TrivS$ получена из $\Triv$  движением $R_3$.
		
Предположим, что перекрестки $x_1$ и $x_2$ у диаграмм $\Triv$ и $\B$ совпадают. Тогда выберем такой $x$, при котором $\B=\Triv_{x}$ и в $x$ нет перекрестков  $x_1$ и $x_2$. По утверждению \ref{p:hardest}   $0=\c(\Triv)=\c(\TrivS)=\c(\Triv_{R_3})$. Тогда по утверждению \ref{p:Mov}  $\c(\Triv_x)=\c(\Triv_{R_3,x'})=\c(\B)=\c(\B_{R_3})$. 

Предположим, что один перекресток $x_i$ из $x_1$ и $x_2$ диаграммы $\Triv$ не совпадает с перекрестком диаграммы $\B$. Нетрудно видеть, что неупорядоченные диаграммы $\Triv_{[x_i]}$ и $\TrivS_{[x_i']}$ эквивалентны. Имеем
$$\c(\Triv_{x_i})=\c(\Triv)+\s(x_i,\Triv_{x_i})\b(\Triv_{[x_i]})=\c(\TrivS)+\s(x_i',\TrivS_{x_i'})\b(\TrivS_{[x_i']})=\c(\TrivS_{x_i'}).$$
Выберем такое $x$, при котором $\B=\Triv_{x_i,x}$ и в $x$ нет перекрестков $x_1$ и $x_2$. Имеем $\c(\Triv_{x_i})=\c(\TrivS_{x_i'})=\c(\Triv_{x_i,R_3})$. Тогда по утверждению \ref{p:Mov}  $\c(\Triv_{x_i,x})=\c(\Triv_{x_i,R_3,x'})=\c(\B)=\c(\B_{R_3})$.

Предположим, что перекрестки $x_1$ и $x_2$ у диаграмм $\Triv$ и $\B$ не совпадают. Нетрудно видеть, что для какого-то $i$ диаграммы $\Triv_{[x_i]}$ и $\TrivS_{[x_i']}$ эквивалентны. Имеем
$$\c(\Triv_{x_i})=\c(\Triv)+\s(x_i,\Triv_{x_i})\b(\Triv_{[x_i]})=\c(\TrivS)+\s(x_i',\TrivS_{x_i'})\b(\TrivS_{[x_i']})=\c(\TrivS_{x_i'}).$$
Нетрудно видеть, что для $j\neq i$ диаграммы $\Triv_{x_i,[x_j]}$ и $\TrivS_{x_i',[x_j']}$ эквивалентны. Имеем
$$\c(\Triv_{x_i,x_j})=\c(\Triv_{x_i})+\s(x_j,\Triv_{x_i,x_j})\b(\Triv_{x_i,[x_j]})=\c(\TrivS_{x_i'})+\s(x_j',\TrivS_{x_i',x_j'})\b(\TrivS_{x_i',[x_j']})=\c(\TrivS_{x_i',x_j'}).$$
Выберем такое $x$, при котором $\B=\Triv_{x_i,x_j,x}$ и в $x$ нет перекрестков $x_1$ и $x_2$. Имеем $\c(\Triv_{x_i,x_j})=\c(\TrivS_{x_i',x_j'})=\c(\Triv_{x_i,x_j,R_3})$. Тогда по утверждению \ref{p:Mov}  $\c(\Triv_{x_i,x_j,x})=\c(\Triv_{x_i,x_j,R_3,x'})=\c(\B_{R_3})=\c(\B)$. $\Box$

\end{document}